\documentclass[12pt]{amsart}
\usepackage{amsmath, amssymb, amsfonts, amsthm,amscd}
\usepackage{mathrsfs,bbm}
\usepackage{txfonts}
\usepackage[all]{xy}
\usepackage{nicefrac}
\newtheorem{theorem}{Theorem}

\newtheorem{corollary}[theorem]{Corollary}

\newtheorem{proposition}[theorem]{Proposition}

\newtheorem{lemma}[theorem]{Lemma}

\theoremstyle{definition}
\newtheorem{definition}[theorem]{Definition}
\newtheorem{remark}[theorem]{Remark}

\newcommand\Gr{\operatorname{Gr}}
\newcommand\SL{\operatorname{SL}}

\newcommand\res{\operatorname{res}}

\newcommand\Ind{\operatorname{Ind}}

  \newcommand{\fg}{\mathfrak{g}}
\newcommand{\ft}{\mathfrak{t}}

\newcommand{\bc}{\mathbb{C}}

\newcommand{\g}{\mathfrak{g}}
\renewcommand{\ni}{\noindent}
 \newcommand{\C}{\mathbb{C}}

\def\fg{{\mathfrak g}}

\def\ft{{\mathfrak t}}
\def\Square{}
\def\svert{\vert}
\def\coloneq{:=}

\title[Category ${\mathscr{C}}_{k}$ of multi-loop algebra representations]
{ Category ${\mathscr{C}}_{k}$ of multi-loop algebra representations versus modular representations: Questions of G. Lusztig}

\author{SHRAWAN KUMAR}
\date{}

\address{S. Kumar: Department of Mathematics, University of North Carolina, Chapel Hill, NC 27599-3250, USA}
\email{shrawan@email.unc.edu} 

\begin{document}
\maketitle{}

\section{Abstract} Lusztig defined an abelian  category ${\mathscr{C}}_{k}$ of a class of representations of a multi-loop algebra
and asked various questions connecting it to the modular representation theory of simple algebraic groups in char. $p>0$. The aim of this paper is to show that some of these questions have negative answer. 
\vskip4ex

\section{Introduction}
Let $G$ be a connected simply-connected simple algebraic group over an
algebraically closed field of characteristic $p > 0$ with a fixed maximal torus $T$ and
a fixed Borel subgroup $B$ containing $T$. Let $X_+$  be the set of dominant characters of
$T$.  Consider the category $\mathfrak{C}=\mathfrak{C}(G)$  of finite dimensional rational representations of $G$.
The simple objects of $\mathfrak{C}$, up to isomorphism, 
are indexed by $X_+$; let $L_\lambda$ be the simple object indexed by $\lambda \in X_+$. Let $E_\lambda^0$ 
 be
the Weyl module indexed by $\lambda \in X_+$. The Weyl modules form another basis of the
Grothendieck group $\mathcal{G}(\mathfrak{C})$ of $\mathfrak{C}$.  Hence, for any $\lambda \in X_+$, we can write
$$L_\lambda = \sum_{\mu\in X_+} c_{\mu, \lambda}\,  E_\mu^0,$$
where $c_{\mu, \lambda}$ are integers,  and they are zero for all but finitely many $\mu$.
It is of considerable
interest to understand the character of each $L_\lambda$ or, equivalenty, to understand the
integers $c_{\mu, \lambda}$ (since the characters of $E_\mu^0$ are given by the Weyl character formula). 
By a famous conjecture of Lusztig [Lu1], the integers $c_{\mu, \lambda}$ for $\lambda$ in a finite subset 
of $X_+$ containing the restricted weights $X_+^{\res} : = \{\mu\in X_+: \mu(\alpha_i^\vee) < p \,\text{ for all the simple coroots}\,  
\alpha_i^\vee\}$, are given in terms of the Kazhdan-Lusztig polynomials of the affine Weyl group of the Langlands dual of $G$,  assuming that $p$ is sufficiently large relative
to the type of $G$. (By the Steinberg tensor product theorem, this leads to a  formula for $c_{\mu, \lambda}$
for any $\lambda\in X_+$.)
 This conjecture has been proved for $p$ `very
large' by Andersen-Jantzen-Soergel [AJS]. Fiebig \cite{F} proved the conjecture still for very large but explicit bound on $p$. 

Now, more recently, Lusztig [Lu2] has formulated a conjecture in
such a way that the tensor product theorem is not used. Namely, for any $\lambda\in X_+$
and any integer $k\geq 0$,  he defined an element $E^k_\lambda\in \mathcal{G}(\mathfrak{C})$
by induction on $k$. When $k = 0$, $E^0_\lambda$ is already defined above. Define $E^1_\lambda$ 
to be the reduction mod $p$ of the simple module with highest weight $\lambda$ 
of the quantum group associated to $G$ at a $p$-th root of $1$. Then,
$$E_\lambda^1= \sum_{\mu\in X_+}\,\mathcal{P}_{\mu, \lambda} E^0_\mu,\,\,\,\text{for $\mathcal{P}_{\mu, \lambda}\in \mathbb{Z}$},$$
where $\mathcal{P}_{\mu, \lambda}$ is explicitly known in terms of the Kazhdan-Lusztig polynomials of the affine Weyl group due to the works of Kazhdan-Lusztig [KL0-4] and Kashiwara-Tanisaki [KT]. Express $\lambda$ in terms of its $p$-power expansion:
$\lambda = \sum_{r\geq 0}\, p^r\lambda_r,\,\,\,\text{where $\lambda_r\in X_+^{\res}$}.$
Following Lusztig [Lu2], define
$$ E_\lambda^2= \sum_{\mu\in X_+:\mu-\lambda_0\in pX}\,\mathcal{P}_{(\mu-\lambda_0)/p, (\lambda- \lambda_0)/p} E^1_\mu,\,\,\,
 E_\lambda^3= \sum_{\mu\in X_+:\mu-\lambda_0-p\lambda_1\in p^2X}\,\mathcal{P}_{(\mu-\lambda_0-p\lambda_1)/p^2, (\lambda- \lambda_0-p\lambda_1)/p^2} E^2_\mu,$$
where $X$ is the set of all the characters of $T$, and continue this way to define $ E_\lambda^k$ for any $k$. 
Then, $\{E^k_\lambda: \lambda \in X_+\}$ is a $\mathbb{Z}$-basis of $\mathcal{G}(\mathfrak{C})$ for any $k\geq 0$. Moreover, each $E^k_\lambda$ is, in fact, a $G$-module (not just a virtual module). Further, for large $k$, 
$E_\lambda^k = E_\lambda^{k+1}=  E_\lambda^{k+2}= \cdots.$
This common value is denoted by $ E_\lambda^{\infty}$.
If $p$ is `very large' compared to the type of $G$, then  $ E_\lambda^{\infty} = L_\lambda$. Thus, this gives an explicit successive approximation to construct $L_\lambda$ for large $p$ (cf. \cite{Lu2} for all these results). 
\vskip1ex

Now, Lusztig has defined parallel objects in characteristic zero in terms of the representation theory of multi-loop algebras [Lu3]. Specifically, let $\mathfrak{g}$ be the (complex) Lie algebra of a simple group over the complex numbers of the same type as $G$ and let $\langle\, ,\,\rangle$ be the invariant (symmetric) bilinear form on $\mathfrak{g}$ normalized so that the induced  form on the dual $\mathfrak{t}^*$ of the Cartan subalgebra $\mathfrak{t}$ of $\fg$ satisfies: $\langle \theta, \theta\rangle =2$ for the highest root $\theta$. For any $k\geq 0$, define the Laurent polynomial ring 
$$A_k:= \mathbb{C}[t_1^{\pm 1}, \dots, t_k^{\pm 1}]$$ and let 
$$\tilde{\mathfrak{g}}_k := (\mathfrak{g}\otimes A_k) \oplus (\C c_1 \oplus \dots \oplus \C c_k)$$
 with bracket defined by
\begin{align*}[x\otimes {\bf t}^{\bf n}, y\otimes {\bf t}^{\bf m}] = [x, y]\otimes {\bf t}^{{\bf n}+{\bf m}}+ &\delta_{{\bf n}+{\bf m}, 0}\langle x, y\rangle (n_1c_1+ \cdots + n_kc_k),\,\,\text{for $x, y\in \mathfrak{g}$},\\
& {\bf n} = (n_1, \dots, n_k)\in \mathbb{Z}^k\,\,\text{ and}\,\, {\bf m} = (m_1, \dots, m_k)\in \mathbb{Z}^k,
\end{align*}
where ${\bf t}^{\bf n} := t_1^{n_1}\cdots t_k^{n_k}$ and $c_i$ are central elements. 

Let $\mbox{}^{\prime}{\mathscr{C}}_{k}$ be the category of $\tilde{\mathfrak{g}}_k$-modules on which $c_i$ acts by $-p^i -h^\vee$ for $1\leq i\leq k$, where $h^\vee$ is the dual Coxeter number of $\fg$. (In particular, $\mbox{}^{\prime}{\mathscr{C}}_{0}$ is the category of ${\mathfrak{g}}$-modules.) 
For any $M\in \mbox{}^{\prime}{\mathscr{C}}_{k-1}$ (with $k \geq 1$) define a $\tilde{\mathfrak{g}}_k$-module $\Ind^k_{k-1}(M)$ as follows:
Let  
$$A^+_k:= \mathbb{C}[t_1^{\pm 1}, \dots, t_{k-1}^{\pm 1}, t_k^1] \subset A_k$$
 and regard $M$ as a $\tilde{\g}_k^+$ -module by letting $\fg \otimes (t_kA_k^+)$ act by zero and $c_i$ acts by $-p^i-h^\vee$, where
 $$\tilde{\g}_k^+:= (\fg\otimes A_k^+)\oplus (\C c_1 \oplus \dots \oplus \C c_k).$$
  Now, let
$$\Ind^k_{k-1}(M) := \Ind^{U(\tilde{\mathfrak{g}}_k)}_{U(\tilde{\g}_k^+)}\,(M).$$
It clearly belongs to the category  $\mbox{}^{\prime}{\mathscr{C}}_{k}$. 

Define the sequence of objects ${\mathscr{E}}^{0}_{\lambda,k}, {\mathscr{E}}^{1}_{\lambda,k},\dots, \mathscr{E}^{k}_{\lambda,k}$ in $\mbox{}^{\prime}{\mathscr{C}}_{k}  $ by induction on $k$ by setting ${\mathscr{E}}^{0}_{\lambda,0}$
as the Weyl (which is the same as irreducible) module $V(\lambda) $ of $\mathfrak{g}$ with highest weight $\lambda$. For $k \geq 1$ and $k'\in \{1, 2, \dots, k-1\}$, set 
$$\mathscr{E}^{k'}_{\lambda,k} := \Ind^k_{k-1}(\mathscr{E}^{k'}_{\lambda,k-1}).$$ 

Lusztig asked the following questions $(Q_1) - (Q_4)$ (cf. [Lu3]). Actually, he termed these questions as his `expectations'. 
\vskip1ex

\ni
{\bf Question $Q_1$ (Lu3):} {\it The $\tilde{\mathfrak{g}}_k$-module  ${\mathscr{E}}^{0}_{\lambda,k}$ has a unique irreducible quotient. Denote it by ${\mathscr{E}}^{k}_{\lambda,k}$. 
(For $k=1$ this is proved in [KL1]. )}

\vskip1ex
Assuming the validity of the above question $(Q_1)$,
let ${\mathscr{C}}_{k}$ be the abelian subcategory of $\mbox{}^{\prime}{\mathscr{C}}_{k}$ consisting of those 
$\tilde{\mathfrak{g}}_k$-modules which have finite length and have all their composition factors in the $\tilde{\mathfrak{g}}_k$-
irreducible modules $\{{\mathscr{E}}^{k}_{\lambda,k}: \lambda\in X_+\}$.

\vskip1ex
\ni
{\bf Question $Q_2$ (Lu3):} {\it Each of the modules ${\mathscr{E}}^{k'}_{\lambda,k}$, where $0\leq k'\leq k$ and $\lambda\in X_+$, lie in ${\mathscr{C}}_{k}$. (For $k=1$ this is proved in [KL1]. )
}
\vskip1ex

\vskip1ex
\ni
{\bf Question $Q_3$ (Lu3):} {\it For $0\leq k'\leq k$ and $\lambda\in X_+$, the matrix expressing 
 ${\mathscr{E}}^{k'}_{\lambda,k}$ in terms of the irreducible objects ${\mathscr{E}}^{k}_{\mu,k}\,(\mu\in X_+)$ in the Grothendieck group $\mathcal{G}({\mathscr{C}}_{k})$ is the same as the matrix expressing $(E_\lambda^{k'})$ in terms of $(E_\mu^{k})$.}
\vskip1ex

\vskip1ex
\ni
{\bf Question $Q_4$ (Lu3):} {\it The category ${\mathscr{C}}_{k}$ is a rigid braided monoidal category. In particular, 
the category ${\mathscr{C}}_{k}$ is equipped with a duality operator.
(For $k=1$ this is proved in [KL1, KL2].) }
\vskip1ex

We enlarge the Lie algebra $\tilde{\fg}_k$ (resp. $\tilde{\fg}_k^+$) by adding the standard derivations $d_1, \dots, d_k$ to get the Lie algebra $\hat{\fg}_k$ (resp. $\hat{\fg}_k^+$) (cf. Definition \ref{defi1.1}). Replacing  $U(\tilde{\mathfrak{g}}_k)$ 
(resp. $U(\tilde{\mathfrak{g}}_k^+)$)
by $U(\hat{\mathfrak{g}}_k)$ (resp. $U(\hat{\mathfrak{g}}_k^+)$) in the definition of 
$\Ind^{U(\tilde{\mathfrak{g}}_k)}_{U(\tilde{\g}_k^+)}\,(M)$, we get the $\hat{\fg}_k$-modules (cf. Definition \ref{defi1.1})
$${\hat{\mathscr{E}}}^{0}_{\lambda,k}, {\hat{\mathscr{E}}}^{1}_{\lambda,k},\dots, \hat{\mathscr{E}}^{k}_{\lambda,k}.$$  

We prove that for any $\lambda\in X_+$, 
the module $\hat{{\mathscr{E}}}^{0}_{\lambda,k}$ (and hence any $\hat{{\mathscr{E}}}^{j}_{\lambda,k}$ for $j<k$) 
has a unique irreducible quotient as a $\hat{\mathfrak{g}}_k$-module  denoted by $\hat{\mathscr{E}}^{k}_{\lambda,k}$ (cf. Lemma \ref{lemma1}). This answers the Question $Q_1$ affirmatively for the enlarged Lie algebra $\hat{\fg}_k$. 

Our main technical result of the paper is the following theorem (cf. Theorem \ref{pps1}
and Remark \ref{remark1}):

\begin{theorem} Take $k\geq 2$. Then, for any non-constant element $\mathcal{Z} \in U\left(\hat{\mathfrak{g}}_k^-\right)$, there exists $p_o>0$ such that  $\left [\mathfrak{g}\otimes t_{k-1}^{r},\mathcal{Z}\right]\neq0$ for all $r \geq p_o$, where 
$\hat{\mathfrak{g}}_k^-$ is defined by the equation \eqref{eqn0.2}.
\end{theorem}

We next prove that each of $\hat{\mathscr{E}}_{\lambda,k}^{k}$ is an irreducible $\tilde{\mathfrak{g}}_{k}$-module for any  $k\geq 0$ (cf. Proposition \ref{pps4}). Even though we do not know if the $\tilde{\mathfrak{g}}_k$-module  ${\mathscr{E}}^{0}_{\lambda,k}$ has a unique irreducible quotient, but it has a `preferred' irreducible $\tilde{\fg}_k$-module  quotient $\hat{\mathscr{E}}_{\lambda,k}^{k}$. Thus, for any $k \geq 0$,
$$ {\mathscr{E}}_{\lambda,k}^{k'} = \hat{\mathscr{E}}_{\lambda,k}^{k'},\,\,\text{for any $\lambda \in X_+$ and $0\leq k'\leq k$}.$$
The following result, which is the main result of the paper, asserts that Question $Q_2$ has a negative answer for any $\fg$ and any $k\geq 2$ (cf. Corollary \ref{cry7} deduced from Theorem \ref{pps6}). Since $Q_2$ has negative answer, Question $Q_3$ does not make sense at least as it is. 

\begin{corollary}
 Lusztig's question $Q_2$  has negative answer for any $\mathfrak{g}$ and any $k\geq2$. Specifically, ${\mathscr{E}}^{k-1}_{\lambda,k} 
 $ (and hence any ${\mathscr{E}}^{k'}_{\lambda,k} 
 $ for $k'<k$) does not belong to the category $\mathscr{C}_k$ unless ${\mathscr{E}}^{k-1}_{\lambda,k}$ is  irreducible as a $\tilde{\mathfrak{g}}_k$-module.
\end{corollary}

In the last section we study the action of the Sugawara operator on any $\hat{\fg}_k$-module $V$ which has its $d_k$-eigenvalues bounded above (cf. Theorem \ref{sugawara}). In fact, our theorem is slightly more general and it extends the corresponding well-known theorem for $k=1$. \vskip3ex

 \noindent
 {\bf Acknowledgements:} I am grateful to G. Lusztig for bringing the above questions to my knowledge and many subsequent conversations.
 
\section{Questions $Q_1$ and $Q_2$ for loop algebras with derivations}
\begin{definition} \label{defi1.1} Let $\fg$ be a simple Lie algebra over $\bc$. Fix $k\geq 1$.
Just as in the definition of affine Kac-Moody Lie algebras one adds a derivation to make the root spaces as well as weight spaces in integrable highest weight modules  finite dimensional, we add derivations $d_1,\dots,d_k$ to $\tilde{\mathfrak{g}}_k$ to get 
\begin{equation} \label{eqn0.1}
\hat{\mathfrak{g}}_k := \tilde{\mathfrak{g}}_k
\oplus\left(\underset{i=1}{\overset{k}{\oplus}}\mathbbm{C}d_i\right),
\end{equation}
where the bracket is defined by (for all $1\leq i, j\leq k)$,
$$
[d_i, x\otimes {\bf t}^{\bf n}] = n_i x\otimes {\bf t}^{{\bf n}},\,\,\text{for $x\in \mathfrak{g}$ and ${\bf n} = (n_1, \dots, n_k)\in \mathbb{Z}^k$},\,\,\text
{and}\,\,
[d_i, d_j] = [d_i, c_j] = 0.$$
We also define
$$
\hat{\mathfrak{g}}^+_k = \tilde{\mathfrak{g}}^+_k
\oplus \left(\underset{i=1}{\overset{k}{\oplus}} \mathbbm{C}d_i\right).
$$
Let $\mbox{}^{\prime}\hat{\mathscr{C}}_{k}$ be the category of $\hat{\mathfrak{g}}_k$-modules on which $c_i$ acts as $-p^i -h^\vee$. For $\lambda\in X_+$,  we define the sequence of objects $\hat{{\mathscr{E}}}^{0}_{\lambda,k}, \hat{{\mathscr{E}}}^{1}_{\lambda,k},\dots, \hat{\mathscr{E}}^{k}_{\lambda,k}$ in $\mbox{}^{\prime}\hat{\mathscr{C}}_{k}  $
exactly by the  definition as in the Introduction replacing $\tilde{\mathfrak{g}}_k$ by $\hat{\mathfrak{g}}_k$. In defining $\hat{\mathscr{E}}_{\lambda,k}^{j}$, we let $d_k$ act trivially on $\hat{\mathscr{E}}_{\lambda,k-1}^{j}$ (for any $0\leq j <k$). ($\hat{\mathscr{E}}_{\lambda,0}^{0}$ is defined to be the Weyl module ${\mathscr{E}}_{\lambda,0}^{0} = V(\lambda)$.)

Let $A_k^-\coloneq t^{-1}_k\mathbbm{C}\left[t_1^{\pm1},\dots,t_{k-1}^{\pm1},t_k^{-1}\right]$ and
\begin{equation}\label{eqn0.2}
\hat{\mathfrak{g}}_k^-=\mathfrak{g}\otimes A_{k}^-.
\end{equation}
\end{definition}

The following lemma answers the Question $Q_1$ affirmatively for $\hat{{\mathscr{E}}}^{0}_{\lambda,k}$ replacing 
${\mathscr{E}}^{0}_{\lambda,k}$.

 \begin {lemma}\label{lemma1} For any $\lambda\in X_+$, 
the module $\hat{{\mathscr{E}}}^{0}_{\lambda,k}$ (and hence any $\hat{{\mathscr{E}}}^{j}_{\lambda,k}$ for $j<k$) 
has a unique irreducible quotient as a $\hat{\mathfrak{g}}_k$-module. Let us denote it by $\hat{\mathscr{E}}^{k}_{\lambda,k}$.
\end {lemma}

\begin{proof}
Let \{$M_{{\alpha}}$\} be the collection of all the proper $\hat{\mathfrak{g}}_k$-submodules of $\hat{\mathscr{E}}_{\lambda,k}^{0}$
and let 
$$
M :=\sum M_{{\alpha}}\subset \hat{\mathscr{E}}^{0}_{\lambda,k}.
$$

Then, $M$ is the maximal proper $\hat{\mathfrak{g}}_k$-submodule of $\hat{\mathscr{E}}^{0}_{\lambda,k}$. Observe that it is a proper submodule since $\hat{\mathscr{E}}^{0}_{\lambda,k}$ is diagonalizable as a $D_k := (\mathbbm{C}d_1 \oplus\dots, \oplus\,\mathbbm{C}d_k)$-module and any non-zero element of the 0-weight space of $\hat{\mathscr{E}}^{0}_{\lambda,k}$ with respect to the action of $D_k$ generates $\hat{\mathscr{E}}^{0}_{\lambda,k}$ as a $\hat{\mathfrak{g}}_k$ -module. Thus, any $M_{\alpha}$ (and hence $M$) cannot contain any non-zero element of the 0-weight space.
\end{proof}

Let $U$ denote the universal enveloping algebra.

\begin{theorem}\label{pps1} Take $k=2$.
For any non-constant element $\mathcal{Z} \in U\left(\hat{\mathfrak{g}}_2^-\right)$, there exists $p_o>0$ (depending upon $\mathcal{Z}$) such that the commutator $\left [\mathfrak{g}\otimes t_1^{r},\mathcal{Z}\right]\neq0$ for all $r \geq p_o$.
\end{theorem}

\begin{proof}
Choose a basis $\left\{x_i\right\}_{i\geq1}$ of $\hat{\mathfrak{g}}^-_2$ of the form
$$
x_i=e_{\beta_{i}}\otimes t_1^{p_i}t_2^{q_i},\,\,\,q_i<0,
$$
where $\left\{e_{\beta}\right\}$ is a basis of $\mathfrak{g}$ consisting of root vectors ($e_\beta$ has root $\beta$) and a basis of the Cartan subalgebra $\ft\subset \mathfrak{g}$ dual to the simple roots $\left\{\alpha_i\right\}$. Write (by renumbering the basis elements) the top homogeneous component $\mathcal{Z}^0$ of $\mathcal{Z}$ in the PBW-basis:
$$
\mathcal{Z}^0 = \sum c_{{\bf d}} x_1^{d_1}\dots x_n^{d_n}, ~ 0\neq c_{{\bf d}}\in \mathbbm{C} \quad \mbox{where}~ {\bf d}=\left(d_1,\dots, d_n\right).
$$
In the above expression we only list those $x_i$ such that $x_i$ appears with exponent $d_i\geq 1$ in at least one monomial ${\bf x}^{{\bf d}}\coloneq x_1^{d_1}\dots x_n^{d_n}$ with $c_{{\bf d}}\neq 0$.
We also index $x_i$ so that
$$
x_1,\dots,x_m \in\quad  \ft\otimes A_2^-\qquad (m\geq 0)
$$
and $x_{m+1},\dots,x_n \in$ {Root spaces} $\otimes A_2^-$.
\vskip1ex

Let the homogeneous component $\mathcal{Z}^o$ be of degree $d_{\mathcal{Z}}$.
\vskip1ex

$\underline{\mbox{Case I}}:$ $m< n$, i.e., $\mathcal{Z}^o\notin U\left( \ft\otimes A_2^-\right)$. 
\vskip1ex

Of course, $\mathcal{Z}^o=\sum_{\svert {\bf d}\svert=d_{\mathcal{Z}}}c_{{\bf d}}x_1^{d_1}\dots x_n^{d_n}$, where $\svert {\bf d}\svert :=d_1+d_2+\cdots+d_n$.
Let $\bar{d}_{m+1}\coloneq \max  \left\{d_{m+1}:c_{{\bf d}}\neq 0\right\}$.

Fix $\left( \bar{d}_1,\dots, \bar{d}_m\right)$ such that $c_{(\bar{d}_1,\dots,\bar{d}_m, \bar{d}_{m+1}, \dots)} \neq 0$ for some $d_{m+2},\dots,d_n$.

Set
\begin{eqnarray*}
&&\bar{d}_{m+2}:= \max \left\{ d_{m+2}:  c_{( \bar{d}_1, \bar{d}_2, \dots, \bar{d}_m,\bar{d}_{m+1},d_{m+2},\dots)}\neq0\mbox{ for some }d_{m+3},\dots,d_n\right\}\\
&&\vdots\\
&&\bar{d}_n  := \max\left\{ d_n : c_{( \bar{d}_1,\dots , \bar{d}_{n-1},d_n)} \neq 0 \right\} .
\end{eqnarray*}
Then, of course, $\bar{d}_n= d_{\mathcal{Z}}-\left( \bar{d}_1+ \cdots + \bar{d}_{n-1}\right)$.
Let $m+1 \leq k \leq n$ be the maximum integer such that
$$
\bar{d}_k\neq 0.
$$
Let $x_k=e_{\beta_k}\otimes t_1^{p_k}t_2^{q_k}$ for a root $\beta_k\neq 0. $
Take $h_o\in \ft$ such that $\beta_k(h_o)\neq 0$. 
Then, we claim that the coefficient of
\begin{equation}\label{eqn5.1}
A:= x_1^{\bar{d}_1}x_2^{\bar{d}_2}\dots x_{k-1}^{\bar{d}_{k-1}} x_{k}^{\bar{d}_{k}-1} x_{k+1}^{0}\dots x_{n}^{0}\cdot \left( e_{\beta_k}\otimes \left ( t_1^{p_k+r}t_2^{q_k}\right)\right )\,\,\text{
in}\,\, \left[ h_o \otimes t_1^r, \mathcal{Z}^o\right]_{d_{\mathcal{Z}}}
\end{equation} 
is nonzero if we take $r>>0$, 
$\left[ h_o \otimes t_1^r, \mathcal{Z}^o\right]_{d_{\mathcal{Z}}}$ denotes the $d_{\mathcal{Z}}$-graded component in $\Gr \left( U\left( \mathfrak{g}\otimes A_2^-\right)\right)$ under the  standard filtration of the enveloping algebra.

To prove the above claim, observe that the coefficient of $A$ in $\left[ h_o \otimes t_1^r, \mathcal{Z}^o\right]_{d_\mathcal{Z}}$ can only come from (by using the definition of $\bar{d}_i$) the commutator of $h_o \otimes t_1^r$ with the monomial $x_1^{\bar{d}_1}x_2^{\bar{d}_2}\dots x_{k-1}^{\bar{d}_{k-1}} x_k^{\bar{d}_k}x_{k+1}^{0}\dots x_n^{0}$ or the monomials: 
$$\left\{ B_j:= x_1^{\bar{d}_1}x_2^{\bar{d}_2}\dots x_{k-1}^{\bar{d}_{k-1}} x_{k}^{\bar{d}_{k}-1}x_{k+1}^{0}\dots x_j^{1}\dots x_n^0\right\}_{k+1 \leq j\leq n}.
$$
Now, the component of $A$ in 
$$
\left[ h_o\otimes t_1^r,x_1^{\bar{d}_1}\dots x_k^{\bar{d}_k} x_{k+1}^{0}\dots x_n^0\right]_{d_{\mathcal{Z}}}=\beta_k \left(h_o\right)\cdot A .
$$
Further, for any $k+1\leq j\leq n,$
let $x_j=e_{\beta_j}\left(t_1^{p_j}t_2^{q_j}\right)$.
Then, if $\left( p_j,q_j\right)\neq \left( p_k,q_k\right)$, then 
the component of $A$ in $\left[h_o \otimes t_1^r,B_j\right]_{d_{\mathcal{Z}}}$ is clearly zero. But, if $\left( p_j,q_j\right)=\left( p_k,q_k\right)$, then the root $\beta_j$ will have to be different from $\beta_k$. 
Thus, again the component of $A$ in 
$\left[h_o \otimes t_1^r,B_j\right]_{d_{\mathcal{Z}}}$ is zero. 
Hence, the coefficient of $A$ in $\left[h_o \otimes t_1^r,\mathcal{Z}^o\right]_{d_{\mathcal{Z}}}$ is nonzero, proving the claim \eqref{eqn5.1}. Thus, we have proved the  theorem in this Case I by showing that 
$\left[h_o \otimes t_1^r,\mathcal{Z}^o\right]_{d_{\mathcal{Z}}}\neq 0$ which of course implies that $\left[h_o \otimes t_1^r,\mathcal{Z}\right]\neq 0$. 
\vskip1ex

$\underline{\mbox{Case II}}$: $m=n$, i.e., $\mathcal{Z}^o\in U\left( \ft\otimes A_2^-\right)$. 
\vskip1ex

Again we write
$$
\mathcal{Z}^o=\sum c_{{\bf d}} x_1^{d_1}\dots x_m^{d_m} , 
$$
where we now have each $x_i\in \ft\otimes A_2^-$. 
Similar to the proof in Case I, define 
\begin{eqnarray*}
\bar{d}_1&=& 
\max \left\{ d_1 : c_{{\bf d}}\neq 0\right\}\\
\bar{d}_2&=& \max \left\{ d_2 : c_{(\bar{d}_1,{d}_2,\dots)}\neq 0\quad \mbox{for some }d_3,\dots,d_
m \right\},\\
\vdots&&\\
\bar{d}_m&=& \max \left\{ d_m : c_{(\bar{d}_1,\dots, \bar{d}_{m-1},d_m)}\neq 0\right\}. 
\end{eqnarray*}
Of course, $\bar{d}_m=d_{\mathcal{Z}}-\left( \bar{d}_1+\bar{d}_2+\cdots +\bar{d}_{m-1}\right)$. 
Let $1\leq k \leq m$ be the maximum integer such that $\bar{d}_k\neq 0$. 
Let $x_k=h_k \otimes t _1^{p_k}t_2^{q_k}$. 
Take the simple root $\beta_k$ such that $\beta_k(h_k)\neq 0$. (Observe that we have taken the basis of $\ft$ dual to the simple roots.) Then, we claim that the coefficient of 
\begin{equation}\label{eqn5.2}
A' := x_1^{\bar{d}_1}\dots x_{k-1}^{\bar{d}_{k-1}}x_k^{\bar{d}_{k}-1}x^0_{k+1}\dots x^0_m \left( e_{\beta_k } \otimes t _1^{r+p_k}t_ 2^{q_k}\right) \,\,\text{
in}\,\, \left[ e_{\beta_k}\otimes t_1^r, \mathcal{Z}^o\right]_{d_{\mathcal{Z}}}
\end{equation}
 is nonzero for $r>>0$.

Again by a proof similar to that of Case I, the coefficient of $A'$ in $\left[ e_{\beta_k}\otimes t_1^r, \mathcal{Z}^o\right]_{d_{\mathcal{Z}}}$ can only come from the commutator of $e_{\beta_k}\otimes t_1^r$ with the monomial $x_1^{\bar{d}_1}\dots x_k^{\bar{d}_k}x_{k+1}^0\dots x_m^0$ or the monomials
$$
\left\{ B_j' := x_1^{\bar{d}_1} \dots x_{k-1}^{\bar{d}_{k-1}} x_{k}^{\bar{d}_{k}-1} x_{k+1}^0 \dots x_j^1 \dots x_m^0 \right\} _{k+1\leq j\leq m}. 
$$
Now, the component of $A'$ in 
$\left[ e_{\beta_k}\otimes t _1^r ,x_1^{\bar{d}_1} \dots x_{m}^{\bar{d}_m}\right]_{d_{\mathcal{Z}}}$ is easily seen to be equal to $-A'$.
Further, for any $k+1\leq j\leq m$, let 
$$
x_j=h_j \otimes t _1^{p_j}t _2^{q_j}. 
$$
If $( p_j,q_j)\neq ( p_k, q_k)$, then, of course,  the component of $A'$ in $\left[ e_{\beta_k}\otimes t_1^r, B_j'\right]_{d_{\mathcal{Z}}}$ is clearly zero. But, if  $( p_j,q_j)= ( p_k, q_k)$, then $h_j\neq h_k$ and hence $\beta_k(h_j)=0$. 
Thus, again the coefficient of $A'$ in $\left[ e_{\beta_k}\otimes t_1^r, B_j'\right]_{d_{\mathcal{Z}}}$ is zero. 
Hence, the coefficient of $A'$ in $\left[ e_{\beta_k}\otimes t_1^r, \mathcal{Z}^o\right]_{d_{\mathcal{Z}}}$ is nonzero, proving the claim \eqref{eqn5.2}. Thus, we have proved the theorem in this case as well. This completes the proof of the theorem. $\Square$
\end{proof}
\begin{remark} \label{remark1} The above proof can easily be adapted to prove the following generalization of Theorem \ref{pps1}:

{\it Take $k\geq 2$. Then, for any non-constant element $\mathcal{Z} \in U\left(\hat{\mathfrak{g}}_k^-\right)$, there exists $p_o>0$ such that  $\left [\mathfrak{g}\otimes t_{k-1}^{r},\mathcal{Z}\right]\neq0$ for all $r \geq p_o$.}
\end{remark}

Let $\lambda$ be a dominant integral weight of $\mathfrak{g}$ (i.e., $\lambda\in X_+$)  and let $V(\lambda)$ be the corresponding Weyl ($=$irreducible) module of $\mathfrak{g}$ with highest weight $\lambda$. Realize $V(\lambda)$ as a $\hat{\mathfrak{g}}_1^+$-module by letting $\mathfrak{g}\otimes t_1^d$ act trivially for any $d>0$ and $c_1$ to act by $-p-h^\vee$ (where $h^\vee$ is the dual Coxeter number of $\fg$). We let $d_1$ to act trivially on $V(\lambda)$. 
Define
$$
\hat{M}(\lambda):= \mbox{Ind}^{U(\hat{\mathfrak{g}}_1)}_{U(\hat{\mathfrak{g}}_1^+)}\left( V(\lambda)\right)
$$
and 
$\hat{L}(\lambda)$ its $\hat{\mathfrak{g}}_1$-module (unique) irreducible quotient (cf. Lemma \ref{lemma1}). Observe that for any $v_o\in \hat{M}(\lambda)$ (and hence for $\hat{L}(\lambda)$) there exists $p_o>0$ (depending upon $v_o$) such that 
$$
(\mathfrak{g}\otimes t_1^r) \cdot v_o =0 \quad \mbox{for all }\quad r\geq p_o.
$$

As a corollary of Theorem \ref{pps1}, we get the following. 

\begin{corollary}\label{cor2} Let $\hat{L}(\lambda)$ be an irreducible $\hat{\mathfrak{g}}_1$-module as above. Then, for any proper $\hat{\mathfrak{g}}_2$-submodule
$$
M \underset{\neq}{\subset}
 \mbox{Ind}^{U(\hat{\mathfrak{g}}_2)}_{U(\hat{\mathfrak{g}}_2^+)}\left( \hat{L}(\lambda)\right)
$$
and for any nonzero vector $v\in M$, there exists $p_o=p_o(v)>0$ such that 
$$
\left( \mathfrak{g}\otimes t_1^r\right) \cdot v \neq 0\quad \mbox{for all } r\geq p_o.
$$
\end{corollary}

\begin{proof}
 Write 
$v=\sum^N_{i=1} z_i\otimes \omega_i$, for $z_i\in U\left( \hat{\mathfrak{g}}^-_2\right)$ and some linearly independent elements $\omega_i\in \hat{L}(\lambda)$. 
Then, for any $x\in \mathfrak{g}$, 
\begin{eqnarray*}
\left( x \otimes t_1^r\right)\cdot v &=& \sum_i \left( x \otimes t_1^r\right ) z_i \otimes \omega _i\\
&=& \sum_i \left [x\otimes t _1^r , z_i\right] \otimes \omega_i
+ \sum_i z_i \otimes \left( x \otimes t _1^r \right) \cdot \omega _i\\
&=& \sum_i \left[ x\otimes t _1^r , z_i\right] \otimes \omega_i,\quad \mbox{for } r>>0.
\end{eqnarray*}
But,  $z_i$ is not a scalar for some $i$ since $M$ is a proper submodule of $\mbox{Ind}^{U(\hat{\mathfrak{g}}_2)}_{U(\hat{\mathfrak{g}}_2^+)}\left( \hat{L}(\lambda)\right)$ and $\hat{L}(\lambda)$ is an irreducible $\hat{\mathfrak{g}}_1$-module. 
Thus, the corollary follows from Theorem \ref{pps1}.  $\Square$
\end{proof}

Let $\hat{\mathscr{E}}^2_{\lambda,2}$ be the $\hat{\mathfrak{g}}_2$-module irreducible quotient of $\mbox{Ind}^{U(\hat{\mathfrak{g}}_2)}_{U(\hat{\mathfrak{g}}_2^+)}\left( \hat{L}(\lambda)\right)$ (cf. Lemma \ref{lemma1}). For any integers $m,n$, let $\hat{\mathscr{E}}^2_{\lambda,2}(m,n)$ be the same module as $\hat{\mathscr{E}}^2_{\lambda,2}$ except that we shift the $d_1,d_2$ weights of the highest weight vector $v_{\lambda}\in \hat{L}(\lambda)$ to $(m,n)$ respectively and extend the $d_1,d_2$ action compatibly to $\hat{\mathscr{E}}^2_{\lambda,2}$. 
\vskip1ex

As a corollary of Corollary \ref{cor2}, we get the following. 

\begin{corollary}\label{cor3} Let $\hat{L}(\lambda)$ be an irreducible $\hat{\mathfrak{g}}_1$-module as above (for $\lambda \in X_+$). Then, if $M\subset \mbox{Ind}^{U(\hat{\mathfrak{g}}_2)}_{U(\hat{\mathfrak{g}}_2^+)}\left(\hat{ L}(\lambda)\right)$ is an irreducible proper $\hat{\mathfrak{g}}_2$-submodule, then $M$ cannot be isomorphic as a $\hat{\mathfrak{g}}_2$-module with the irreducible $\hat{\mathfrak{g}}_2$-module $\hat{\mathscr{E}}^2_{\mu,2}(m,n)$ for any  $\mu\in X_+$ and any pair of integers $(m,n)$.

Thus, if $\mbox{Ind}^{U(\hat{\mathfrak{g}}_2)}_{U(\hat{\mathfrak{g}}_2^+)}\left( \hat{L}(\lambda)\right)$ has finite length but not irreducible; in particular, it has an irreducible proper $\hat{\mathfrak{g}}_2$-submodule $M$, then $M$ {\bf cannot} be  $\hat{\mathfrak{g}}_2$-module isomorphic with $\hat{\mathscr{E}}^2_{\mu,2}(m,n)$ for any $\mu\in X_+$ and any pair of integers $(m,n)$. 
\end{corollary}

\begin{proof} Assume, if possible, that there is a $\hat{\mathfrak{g}}_2$-module isomorphism
$$
\Phi: M \rightarrow \hat{\mathscr{E}}^2_{\mu,2}(m,n)\quad \mbox{for some } (m,n)\in \mathbbm{Z}^2. 
$$

Let $M^o\subset M$ be the eigenspace of $d_2$ with the maximum $d_2$-eigenvalue say $\alpha_o\leq 0$. (Observe that the $d_2$-eigenvalues of $\mbox{Ind}^{U(\hat{\mathfrak{g}}_2)}_{U(\hat{\mathfrak{g}}_2^+)}\left( \hat{L}(\lambda)\right)$ are $\leq 0$.)
Let
$
m_{\mu}=\Phi^{-1}(v_{\mu})$, where $v_\mu$ is the highest weight vector of $\hat{L}(\mu)$. Then, $m_\mu \in M^o$. 
Now, we assign $d_1,d_2$-eigenvalues of $v_\mu$ to be that of $m_{\mu}$. 
Clearly, the $d_2$-eigenspace of $\hat{\mathscr{E}}^2_{\mu,2}(m,n)$ with maximum eigenvalue is $\hat{L}(\mu)$. Thus, the isomorphism $\Phi$ would induce an isomorphism 
$$\Phi^o: M^o \overset{\sim}{\rightarrow} \hat{L}(\mu)\,\,\text{ as $\hat{\mathfrak{g}}_2^+$-modules}. 
$$
By Corollary \ref{cor2}, any nonzero vector ${v}^o\in M^o$ satisfies: 
\begin{equation} \label{eqn5.3} \left( \mathfrak{g}\otimes t _1^r\right ) \cdot v^o \neq 0\,\,\,\text{ for all $r>>0$}.
\end{equation}
But, since $\hat{L}(\mu)$ has $d_1$-eigenvalues bounded above, no vector in $\hat{L}(\mu)$ satisfies \eqref{eqn5.3}. 
This contradicts the existence of $\Phi$, proving the corollary. $\Square$ 
\end{proof}

\begin{remark}
 The above argument can easily be extended to prove similar negative results for $\hat{\mathscr{E}}^{k'}_{\lambda,k}$ for $k\geq 2$ and $0\leq k'<k$ as long as   $\hat{\mathscr{E}}^{k'}_{\lambda,k}$ is not irreducible $\hat{\fg}_k$-module.
 
 \end{remark}

\section{Questions $Q_1$ and $Q_2$ for multi-loop algebras}

\begin{proposition}\label{pps4}
Let $k\geq1$. Assume that $\hat{\mathscr{E}}_{\lambda,k-1}^{k-1}$ is an irreducible $\tilde{\mathfrak{g}}_{k-1}$-module. 
Let $\hat{M}_k\coloneq\hat{M}_k(\lambda)\subset \hat{\mathscr{E}}^0_{\lambda,k}$ be the unique  proper maximal $\hat{\mathfrak{g}}_k$-submodule (cf. Lemma \ref{lemma1}). Then, $\hat{M}_k$ is also maximal as a $\tilde{\mathfrak{g}}_k$-submodule. 

Thus, $\hat{\mathscr{E}}_{\lambda,k}^{k}$ is an irreducible $\tilde{\mathfrak{g}}_{k}$-module as well. 

Hence, by induction, each of $\hat{\mathscr{E}}_{\lambda,k}^{k}$ is an irreducible $\tilde{\mathfrak{g}}_{k}$-module for any  $k\geq 0$. Thus, for any $k\geq 0$, taking the preferred choice for ${\mathscr{E}}_{\lambda,k}^{k}$ to be
$\hat{\mathscr{E}}_{\lambda,k}^{k}$, 
$$ {\mathscr{E}}_{\lambda,k}^{k'} = \hat{\mathscr{E}}_{\lambda,k}^{k'},\,\,\text{for any $\lambda \in X_+$ and $0\leq k'\leq k$}.$$
\end{proposition}

\begin{proof}
Let $\hat{M}_k\subset \tilde{M}_k\underset{\neq}{\subset}\hat{\mathscr{E}}^0_{\lambda,k}$ be a $\tilde{\mathfrak{g}}_{k}$-submodule. Let $\hat{L}_k=\hat{\mathscr{E}}_{\lambda,k}^{k}$ be the $\hat{\mathfrak{g}}_{k}$-module irreducible quotient $\hat{\mathscr{E}}_{\lambda,k}^{0}/ \hat{M}_k$. Then,
$
\tilde{N}_k \coloneq \frac{\tilde{M}_k}{\hat{M}_k} \underset{\neq}{\subset} \hat{L}_k$ is a $\tilde{\mathfrak{g}}_k$-submodule. 
Assume, if possible, that $\tilde{N}_k$ is nonzero. Observe that $d_k$ acts on $\hat{L}_k$ with eigenvalues bounded above by 0. For nonzero $v\in \tilde{N}_k$, decompose it as a sum of eigenvectors under $d_k$:
\begin{equation}\label{eqn2.1}
v=\sum_{i\geq0}v_i,\,\,\text{where } \,\,d_kv_i=-iv_i. 
\end{equation}
(Since $\hat{L}_k$ is a $\hat{\mathfrak{g}}_k$-module, each $v_i\in \hat{L}_k$.)
Now, define
\begin{equation} \label{eqn2.2} \svert v\svert _k=\left\{ \Sigma i: v_i \neq 0\right\}.\end{equation}
Choose a nonzero $\overset{o}{v}\in \tilde{N}_k$ such that $\svert \overset{o}{v} \svert_k\leq \svert v\svert_k$, for all nonzero $v\in \tilde{N}_k$. 
If $\svert \overset{o}{v} \svert_k>0$, then at least one of the components $ \overset{o}{v}_{i_o}\neq0$ with $i_o>0$. 
For any $x\in \mathfrak{g}$ and ${\bf p}=\left(p_1,\dots,p_{k}\right)\in\mathbbm{Z}^{k-1}\times\mathbbm{Z}_{>0}$, $\left(x\otimes {\bf t}^{{\bf p}}\right)\cdot \overset{o}{v}=0$,
because of the minimality of $\svert \overset{o}{v}\svert_{k}$ in $\tilde{N}_{k}$. 
In particular,
$$
\left(x\otimes {\bf t}^{{\bf p}}\right)\cdot \overset{o}{v}_{i_o}=0.
$$
Thus, the $\hat{\mathfrak{g}}_k$-submodule $\hat{Q}_k\subset \hat{L}_{k}$ generated by $\overset{o}{v}_{i_o}$ satisfies:
$$
\svert v\svert_k\geq \svert \overset{o}{v}_{{i}_o}\svert_k=i_o>0, \mbox{ for all } v\in \hat{Q}_k. 
$$
But, clearly $\hat{L}_k$ contains nonzero elements $v'$ with $\svert v'\svert_k=0$ coming from the image of $1\otimes\hat{\mathscr{E}}^{0}_{\lambda,k-1}$ under the projection $\hat{\mathscr{E}}^{0}_{\lambda,k}\rightarrow\hat{L}_k$. 
Hence, the $\hat{\mathfrak{g}}_k$-submodule 
$$
0\neq \hat{Q}_k \underset{\neq}{\subset}\hat{L}_k.
$$
This contradicts the $\hat{\mathfrak{g}}_k$-module irreducibility of $\hat{L}_k$ and hence
$\tilde{N}_k=0$ in the case $\svert \overset{o}{v}\svert_k>0$. 

So, assume now that $\svert \overset{o}{v}\svert_k=0$. Thus, $\overset{o}{v}$ itself is an $d_k$-eigenvetor of eigenvalue 0. But, the $d_k$-eigenspace of eigenvalue 0 in $\hat{L}_k$ clearly is equal to the image of $1\otimes \hat{\mathscr{E}}^0_{\lambda,k-1}$ under the projection
$
\hat{\mathscr{E}}^0_{\lambda,k}\rightarrow \hat{L}_k. 
$
In fact, $\hat{L}_k$ being the unique $\hat{\mathfrak{g}}_k$-module irreducible quotient of $\hat{\mathscr{E}}^0_{\lambda,k}$, it is also an irreducible quotient:
$$
\hat{\mathscr{E}}^{k-1}_{\lambda,k}\twoheadrightarrow \hat{L}_k, 
$$
and $d_k$-eigenspace of eigenvalue 0 in $\hat{L}_k$ is equal to the image of $1\otimes\hat{\mathscr{E}}_{\lambda, k-1}^{k-1}$. Since $\hat{\mathscr{E}}_{\lambda, k-1}^{k-1}$ is an irreducible $\tilde{\mathfrak{g}}_{k-1}$-module by assumption, we get that $\tilde{\mathfrak{g}}_{k-1}$-submodule of $\tilde{N}_k$ generated by $\overset{o}{v}$
must contain the full image $\hat{L}_\mu$ of $1\otimes\hat{\mathscr{E}}_{\lambda, k-1}^{k-1}$ in $\hat{L}_k$. 
But, then
$$
\tilde{N}_k=\hat{L}_k.
$$

This is a contradiction, since $\tilde{N}_k$ was assumed to be properly contained in $\hat{L}_k$. Hence, $\tilde{N}_k=0$ in this case as well. This complete the proof of the proposition. $\Square$
\end{proof}

\begin{lemma}\label{lem5}
Any $\tilde{\mathfrak{g}}_k$-module automorphism $\varphi$ of $\hat{\mathscr{E}}_{\lambda ,k}^{k}$ is the identity automorphism up to a nonzero scalar (for any $k\geq 0$). In particular, it is a $\hat{\mathfrak{g}}_k$-module automorphism. 
\end{lemma}

\begin{proof}
 The lemma is clearly true for $k=0$ (by Schur's Lemma). By induction, assume the validity of the lemma for $k-1$ and take $k\geq1$. Let
$$
M^o=\left\{v\in\hat{\mathscr{E}}^k_{\lambda,k}: \left( \mathfrak{g}\otimes {\bf t}^{{\bf p}}\right)\cdot v=0 \quad \mbox{for all } {\bf p}\in \mathbbm{Z}^{k-1}\times \mathbbm{Z}_{>0}\right\}.
$$
Clearly, $M^o$  is stable under the action of $d_k$ (in fact, $M^o$ is stable under the action of each $d_i$, $1\leq i\leq k$). 
Assume, if possible, that $M^o$  contains a nonzero $d_k$-eignevector $\overset{o}{v}$ of eigenvalue $<0$. Then, the  $\hat{\mathfrak{g}}_k$-submodule $\overset{o}{V}$ of $\hat{\mathscr{E}}^k_{\lambda,k}$ generated by $\overset{o}{v}$ clearly satisfies:
$$
\overset{o}{V}\underset{\neq}{\subset}\hat{\mathscr{E}}^k_{\lambda,k}.
$$
This is a contradiction since $\hat{\mathscr{E}}^k_{\lambda,k}$ is an irreducible $\hat{\mathfrak{g}}_k$-module. 
Thus,
$M^o =$ Image of $1\otimes\hat{\mathscr{E}}_{\lambda, k-1}^{k-1}$ under the projection 
$
\hat{\mathscr{E}}_{\lambda, k}^{k-1} \twoheadrightarrow \hat{\mathscr{E}}_{\lambda, k}^{k}.
$
Hence, the automorphism $\varphi$ of $\hat{\mathscr{E}}_{\lambda, k}^{k}$ restricts to a $\tilde{\mathfrak{g}}_{k-1}$-module automorphism of $1\otimes\hat{\mathscr{E}}_{\lambda, k-1}^{k-1}$ (which is the identity automorphism up to a scalar by induction). This proves the lemma since $\hat{\mathscr{E}}_{\lambda, k}^{k}$ is generated by the image of $1\otimes\hat{\mathscr{E}}_{\lambda, k-1}^{k-1}$ as a $\tilde{\mathfrak{g}}_k$-module. $\Square  $
\end{proof}

\begin{theorem}\label{pps6}
Let $k\geq1$ and let $M\subset \hat{\mathscr{E}}_{\lambda, k}^{k-1}$ be an irreducble $\hat{\mathfrak{g}}_{k}^{(k)}$-submodule, where $\hat{\mathfrak{g}}_{k}^{(k)}\subset \hat{\mathfrak{g}}_{k}$ is the Lie subalgebra $\tilde{\mathfrak{g}}_{k}\oplus \mathbbm{C}d_k$. Then, $M$ is an irreducible $\tilde{\mathfrak{g}}_{k}$-module. 

Moreover, for $k\geq 2$, if $\hat{\fg}_k^{(k)}$-irredicible $M\underset{\neq}{\subset}\hat{\mathscr{E}}_{\lambda, k}^{k-1}$, then $M$ cannot be isomorphic with $\hat{\mathscr{E}}_{\mu, k}^{k}$ as $\tilde{\mathfrak{g}}_{k}$-modules for any $\mu\in X_+$. 
\end{theorem}

\begin{proof}
We first prove that $M$ is an irreducible $\tilde{\mathfrak{g}}_{k}$-module:

Let $d\left(\leq0\right)$ be the largest $d_k$-eignevalue of $d_k$-eigenvectors in $M$. Let us renormalize the $d_k$-eigenvalues by subtracting $d$ so that the largest $d_k$-eigenvalue of $M$ becomes 0. 
Let $N\subset M$ be a nonzero $\tilde{\mathfrak{g}}_{k}$-submodule. Recall the definition of $\svert v\svert_k$ for nonzero vectors from the identity \eqref{eqn2.2}  of the proof of Proposition \ref{pps4}.  Choose nonzero $\overset{o}{v}\in N$ such that 
\begin{equation} \label{eqn3.1} \svert \overset{o}{v}\svert_k\leq \svert v\svert_k\,\,\text{ for all nonzero $v\in N$}. 
\end{equation}

Then, for any $x\in \mathfrak{g}$ and ${\bf p}=\left(p_1,\dots,p_k\right)\in \mathbbm{Z}^{k-1}\times \mathbbm{Z}_{>0}$ ,  by \eqref{eqn3.1},

\begin{equation} \label{eqn3.2} \left( x\otimes {\bf t}^{{\bf p}}\right)\cdot \overset{o}{v}=0.
\end{equation}
If $\svert \overset{o}{v}\svert_k>0$, choose a nonzero $d_k$-eigen component $\overset{o}{v}_{i_0}$ with $\svert \overset{o}{v}_{i_0}\svert_k= i_o>0$. 
Then, by the equation \eqref{eqn3.2}, for any $x\in \mathfrak{g}$ and ${\bf p}\in \mathbbm{Z}^{k-1}\times \mathbbm{Z}_{>0}$, 
$$
 \left( x \otimes {\bf t}^{{\bf p}}\right)\cdot \overset{o}{v}_{i_o}=0.
$$
Thus, the $\hat{\mathfrak{g}}^{(k)}_k$-submodule $M'$ of $M$ generated by $\overset{o}{v}_{i_o}$ is properly contained in $M$. This contradicts the irreducibility of $M$ as a $\hat{\mathfrak{g}}_k^{(k)}$-module. Thus, we conclude that 
$\svert \overset{o}{v}\svert_k=0$. 
In particular, $\overset{o}{v}$ is a $d_k$-eigenvector. Thus, the $\tilde{\mathfrak{g}}_k$-submodule of $N$ generated by $\overset{o}{v}$ is $d_k$-stable. 
This forces
$$
N=M
$$
from the irreducibility of $M$ as a $\hat{\mathfrak{g}}_k^{(k)}$-module. 
This proves that $M$ is irreducible as a $\tilde{\mathfrak{g}}_k$-module as well. 
\vskip1ex

We now prove the second part of the proposition:
\vskip1ex

If possible, assume that there exists a $\tilde{\mathfrak{g}}_k$-module isomorphism (for some $\mu\in X_+$): 
$$
\varphi: M \overset{\sim}{\rightarrow} \hat{\mathscr{E}}^k_{\mu,k}.
$$
Let
$$
M^o:= \left\{ v\in M : \left( x \otimes {\bf t}^{{\bf p}}\right) \cdot v=0\quad \mbox{for all } x\in\mathfrak{g} \mbox{ and } {\bf p}\in \mathbbm{Z}^{k-1}\times \mathbbm{Z}_{>0}\right\},
$$
and similarly define $\left( \hat{\mathscr{E}}^k_{\mu,k}\right)^o$. 
Then, $\varphi$ induces an isomorphism:
$$
M^o\simeq \left( \hat{\mathscr{E}}^k_{\mu,k}\right)^o.
$$
Irreducibility of $M$ as a $\hat{\mathfrak{g}}^{(k)}_{k}$-module forces
$M^o=M_0$, and similarly for $\hat{\mathscr{E}}^k_{\mu,k}$, where $M_0$ is the $d_k$-eigenspace of M of eigenvalue 0. (Observe that $\hat{\mathscr{E}}^k_{\mu,k}$ is irreducible even as a $\tilde{\mathfrak{g}}_k$-module by Proposition \ref{pps4}.)
Thus, we get an isomorphism of $\tilde{\mathfrak{g}}_{k-1}$-modules (by restricting $\varphi$): 
$$
\varphi_0:M_0\simeq 1 \otimes \hat{\mathscr{E}}^{k-1}_{\mu,k-1}.
$$
For any nonzero $v\in M_0$, there exists $p_o=p_o(v)>0$ such that 
\begin{equation} \label{eqn3.3}
\left(\mathfrak{g}\otimes t_{k-1}^r\right)\cdot v \neq 0\,\,\text{ for all $r\geq p_o$}.
\end{equation} 
This is proved for $k=2$ in Corollary \ref{cor2}. (The proof of Corollary \ref{cor2} works equally well if we only assume that  $M$ is a proper $\hat{\fg}_2^{(2)}$-submodule.) The proof for any $k\geq 2$ is the same by using Remark \ref{remark1} and the proof of Corollary \ref{cor2}. However, for any vector $w$ in $\hat{\mathscr{E}}^{k-1}_{\mu,k-1}$ (by the definition) 
\begin{equation} \label{eqn3.4}
 \left(\mathfrak{g}\otimes t_{k-1}^r\right) \cdot w=0\,\,\text{ for all $r>>0$}. 
 \end{equation}
The identities \eqref{eqn3.3} and \eqref{eqn3.4} contradict the existence of the $\tilde{\mathfrak{g}}_{k-1}$-module isomorphism $\varphi_0$. This concludes the proof that $M$ cannot be $\tilde{\mathfrak{g}}_k$-module isomorphic with $\hat{\mathscr{E}}^k_{\mu,k}$ (for any $\mu\in X_+$). $\Square$
\end{proof}

\begin{corollary}\label{cry7}
 Lusztig's question $Q_2$ (cf. Introduction)  has negative answer for any $\mathfrak{g}$ and any $k\geq2$. Specifically, ${\mathscr{E}}^{k-1}_{\lambda,k} 
 $ (and hence any $ {\mathscr{E}}^{k'}_{\lambda,k} 
 $ for $k'<k$) does not belong to the category $\mathscr{C}_k$ unless ${\mathscr{E}}^{k-1}_{\lambda,k}$ is  irreducible as a $\tilde{\mathfrak{g}}_k$-module. (Observe that  $ {\mathscr{E}}^{k'}_{\lambda,k} = \hat{\mathscr{E}}^{k'}_{\lambda,k}$ as in Proposition \ref{pps4}.)  
 \end{corollary}

\begin{proof}
Assume that $\hat{\mathscr{E}}^{k-1}_{\lambda,k}$ has finite length as a $\tilde{\mathfrak{g}}_k$-module; in particular, it would have finite length as  $\hat{\mathfrak{g}}_k^{(k)}$-module. This would insure the existence of an irreducible $\hat{\mathfrak{g}}_k^{(k)}$-submodule $M\subset \hat{\mathscr{E}}^{k-1}_{\lambda,k}$. By Theorem \ref{pps6}, $M$ would be an irreducible $\tilde{\mathfrak{g}}_k$-submodule of $\hat{\mathscr{E}}^{k-1}_{\lambda,k}$; in particular, a $\tilde{\mathfrak{g}}_k$-module composition factor of $\hat{\mathscr{E}}^{k-1}_{\lambda,k}$.  Thus, if $\hat{\mathscr{E}}^{k-1}_{\lambda,k}$ is not irreducible as $\tilde{\fg}_k$-module, then $M 
\underset{\neq}{\subset}\hat{\mathscr{E}}^{k-1}_{\lambda,k}$.
But, then by Theorem \ref{pps6},  $M$ is not $\tilde{\mathfrak{g}}_k$-module isomorphic with $\hat{\mathscr{E}}^{k}_{\mu,k}$ (for any $\mu\in X_+$). This shows that Question $Q_2$ has negative answer. $\Square$
\end{proof}

\section{Sugawara operator for multiloop algebras}

For ${\bf n}=\left( n_1,\dots n_k\right)\in \mathbbm{Z}^k$ also written as  ${\bf n}=n_1\delta_1+\cdots + n_k \delta_k$, denote 
$${\bf n}'={\bf n}-n_k\delta_k .$$
For  $x\in \mathfrak{g}$ and ${\bf n}\in  \mathbbm{Z}^k$, denote
$$
x({\bf n})=x\otimes {\bf t}^{{\bf n}}\in \tilde{\fg}_k.
$$
We abbreviate $d\delta_k$  simply by $d$. Thus, 
 $x({\bf n}+d):= x\otimes ({\bf t}^{{\bf n}} \cdot t_k^d).
 $
 
 Define the {\em Sugawara Operator} 
$$L_0^{(k)}:=\frac{1}{2}\sum_j\, \sum_{d\in \mathbbm{Z}}\, :e_j(-d)\cdot e^j(d):\, \in \hat{U}\left((\fg\otimes t_k^{\pm 1})\oplus \bc c_k\right) \subset \hat{U}(\tilde{\fg}_k),$$
 where $\{e_j\}$ is a basis of $\mathfrak{g}$,   $\left\{e^j\right\}$ is the dual basis with respect to the normalized invariant form $\langle\, ,\,\rangle$ on $\fg$ (as in the Introduction), and $\hat{U}$ is the standard completion along the positive root spaces as in, e.g., \cite[Definition 1.5.8]{K}. 
Even though we have not used the following theorem in rest of the paper, we include it for its future use. This extends the corresponding well-known theorem for $k=1$ (cf. \cite[Proposition 10.1]{KRR}).

\begin{theorem} \label{sugawara} Let $V$ be a  $\hat{\mathfrak{g}}_k^{(k)}$-module  which has $d_k$-eigenvalues bounded above (cf. Theorem \ref{pps6} for the definition of the Lie algebra $\hat{\mathfrak{g}}_k^{(k)}$).  Then, for $x\in\mathfrak{g}$ and ${\bf n}\in \mathbbm{Z}^k$, such that ${\bf n}'\neq {\bf 0}$, $\left[ x({\bf n}),L_0^{(k)}\right]$ as an operator on $V$  is given by : 
\begin{eqnarray*}
\left[x({\bf n}),L_0^{(k)}\right]&=& \sum_j\left(-\sum_{d<\frac{n_k}{2}} e_j \left({\bf n}'+d\right) \cdot \left[x,e^j\right]\left(n_k-d\right)
 + \sum_{d<\frac{n_k}{2}} e_j(d)\cdot \left[ x, e^j\right]\left( {\bf n}-d\right)\right)\\
&&\sum_j\left(-\frac{1}{2} e_j\left( {\bf n}'+ \frac{n_k}{2}\right)\cdot  \left[x,e^j\right]\left(\frac{n_k}{2}\right)+\frac{1}{2}e_j\left(\frac{n_k}{2}\right)\cdot \left[x,e^j\right]\left({\bf n}-\frac{n_k}{2}\right)
\right) + h^\vee n_k x({\bf n}),
\end{eqnarray*}
where we follow the convention that $x\left(\frac{n}{2}\right)=0$ for odd $n$.
\end{theorem}

\begin{proof}
Let $\psi:\mathbbm{R}\rightarrow\{0,1\}$ be the cut-off function defined by 
\begin{eqnarray*}
\psi(t)&=&1\qquad \mbox{for}\quad \svert t\svert \leq 1\\
&=& 0\qquad \mbox{otherwise}.
\end{eqnarray*}
Define the element
$$
L_0^{(k)}(\epsilon)=\frac{1}{2}\sum_j\sum_{d\in\mathbbm{Z}}e_j(-d)\cdot e^j(d)\psi(\epsilon d)\in U\left((\mathfrak{g}\otimes t_k^{\pm1})\oplus \mathbbm{C}c_k\right).
$$
Assume that ${\bf n}'\neq {\bf 0}$ (which is an assumption in the theorem). Then, 
\begin{eqnarray}\label{eqn1.1}
2\left[x({\bf n}),L^{(k)}_0(\epsilon)\right]&=&\sum_j\left( \sum_{\frac{n_k}{2}\leq d}\left[ x, e_j\right] \left( {\bf n}-d\right) \cdot e^j(d) \psi (\epsilon d)
+ \sum_{\frac{n_k}{2}>d} e^j (d)\cdot \left[ x,e_j\right] \left({\bf n}-d\right) \psi \left(\epsilon d\right)\right)\notag\\
&&+\sum_j\left( \sum_{\frac{n_k}{2}>d}\left[\left[ x, e_j\right] \left( {\bf n}-d\right), e^j(d) \right]\psi (\epsilon d)
+ \sum_{-\frac{n_k}{2}\leq d} e_j (-d) \cdot \left[ x,e^j\right] \left({\bf n}+d\right) \psi \left(\epsilon d\right)\right)\notag\\
&&+\sum_j\left(\sum_{-\frac{n_k}{2}> d}\left[ x, e^j\right] \left( {\bf n}+d\right) \cdot e_j(-d) \psi (\epsilon d)
+ \sum_{-\frac{n_k}{2}> d} \left[e_j (-d) ,\left[ x,e^j\right] \left({\bf n}+d\right)\right] \psi \left(\epsilon d\right)\right)\notag\\
&& =\sum_j\left(- \sum_{-\frac{n_k}{2}\leq d} e_j \left( {\bf n}'-d\right) \cdot \left[ x, e^j \right] \left( n_k +d\right) \psi \left(\epsilon \left ( d+n_k\right)\right)\right)\notag\\
&&+\sum_j\left(\sum_{-\frac{n_k}{2}\leq d} e_j \left( -d\right)\cdot  \left[ x, e^j \right] \left( {\bf n} +d\right) \psi \left(\epsilon  d\right)
+\sum_{\frac{n_k}{2}> d} e^j \left( d\right)\cdot \left[ x, e_j \right] \left( {\bf n} -d\right) \psi \left(\epsilon  d\right)\right)\notag\\
&&\sum_j\left(-\sum_{\frac{n_k}{2}> d} e^j \left( {\bf n}'+d\right) \cdot \left[ x, e_j \right] \left( -d +n_k\right) \psi \left(\epsilon  \left ( d-n_k\right)\right)\right)
+2h^\vee x({\bf n})\sum_{-\frac{n_k}{2}\leq d< \frac{n_k}{2}} \psi (\epsilon d), \notag\\
&&\,\,\text{since $\sum_j[[x, e_j], e^j]=2h^\vee x,$ cf. \cite[Lemma 3.3.8]{GW} as $h^\vee :=1+\langle \rho, \theta^\vee\rangle$}\notag\\
&&=\sum_j\left(- \sum_{d\leq \frac{n_k}{2}} e_j \left( {\bf n}' +d\right) \cdot \left[ x, e^j\right] \left( n_k-d \right) \psi \left ( \epsilon \left(-d +n_k\right)\right)\right)\notag\\
&&+\sum_j\left(\sum_{d<\frac{n_k}{2}} e_j (d)\cdot \left[ x,e^j \right]\left( {\bf n}-d \right) \psi \left( \epsilon d\right)
+\sum_{d\leq \frac{n_k}{2}} e_j (d)\cdot \left[ x,e^j \right]\left( {\bf n}-d \right) \psi \left( \epsilon d\right)\right)\notag\\
&&\sum_j\left(- \sum_{d< \frac{n_k}{2}} e_j \left( {\bf n}' +d\right) \cdot \left[ x, e^j\right] \left( -d+n_k \right) \psi \left ( \epsilon \left(d -n_k\right)\right)\right)
+2h^\vee x({\bf n})\sum_{-\frac{n_k}{2}\leq d< \frac{n_k}{2}} \psi (\epsilon d)\notag\\
&&=\sum_j \left(- 2\sum_{d< \frac{n_k}{2}} e_j \left( {\bf n}' +d\right) \cdot \left[ x, e^j\right] \left( n_k-d \right) \psi \left ( \epsilon \left(-d +n_k\right)\right)\right)\\
&&\sum_j\left( 2\sum_{d<\frac{n_k}{2}} e_j (d)\cdot \left[ x,e^j \right]\left( {\bf n}-d \right) \psi \left( \epsilon d\right)
 - e_j\left({\bf n}'+ \frac{n_k}{2}\right) \cdot \left[ x, e^j\right]\left(\frac{n_k}{2}\right) \psi \left(\epsilon \left(\frac{n_k}{2}\right)\right)\right)
 \notag\\
&& \sum_j\left( e_j\left( \frac{n_k}{2}\right) \cdot \left[ x, e^j\right]\left({\bf n}-\frac{n_k}{2}\right) \psi \left(\epsilon \frac{n_k}{2}\right)\right)
+2h^\vee x({\bf n})\sum_{-\frac{n_k}{2}\leq d< \frac{n_k}{2}} \psi (\epsilon d),\notag
\end{eqnarray}
since $\psi$ is symmetric. 
We have used the following relation in the above equation.
$$[x, \sum_j e_j\otimes e^j]= \sum_{j}\left[ x,e_j\right] \otimes e^j + \sum_j e_j \otimes \left[ x, e^j \right]=0\quad \mbox{in}\quad  \mathfrak{g}\otimes \mathfrak{g}.$$
Taking limit as $\epsilon \rightarrow 0$ in the equation \eqref{eqn1.1}, we get 
\begin{eqnarray}\label{eqn1.2}
\left[x({\bf n}), \,\varinjlim_{\epsilon \to 0}
L_0^{(k)}(\epsilon)\right]&=&\sum_j\left(-\sum_{d< \frac{n_k}{2}} e_j \left( {\bf n}'+d\right) \cdot \left[ x, e^j\right] \left(n_k-d\right)
+ \sum_{d< \frac{n_k}{2}} e_j (d)\cdot \left[ x, e^j\right] \left({\bf n}-d\right)\right)\notag\\
&&\sum_j\left(-\frac{1}{2}e_j \left({\bf n}'+ \frac{n_k}{2}\right)\cdot  \left[ x, e^j\right]\left(\frac{n_k}{2}\right)
+\frac{1}{2}e_j \left( \frac{n_k}{2}\right) \cdot \left[ x, e^j\right]\left({\bf n}-\frac{n_k}{2}\right)\right)\notag\\
&&+h^\vee n_kx({\bf n}),
\end{eqnarray}
since, in all the above summands, the $d_k$-eigenvalues of the first term is not more than the $d_k$-eigenvalues of the second term. 
Now, 
\begin{eqnarray}\label{eqn1.3}
2L_0^{(k)}(\epsilon)&=&\sum_j\sum_{d\in{\mathbbm{Z}}} e_j (-d)\cdot  e^j (d) \psi (\epsilon d)\notag\\
&=&\sum_j\left(\sum_{d\geq 0} e_j (-d) \cdot e^j (d) \psi (\epsilon d)
+\sum_{d< 0} e^j (d) \cdot e_j (-d) \psi (\epsilon d)\right)\notag\\
&&\sum_j\left(\sum_{d<0} \left(\left[ e_j,  e^j\right] -dc_k \right)\psi (\epsilon d)\right)\notag\\
&=&\sum_j\left(\sum_{d\geq 0}e_j (-d)\cdot e^j (d)\psi (\epsilon d)
 + \sum_{d<0} e^j (d) \cdot e_j (-d)\psi (\epsilon d)\right)\\
&& - \sum_j\left(\sum _{d<0} d c_k \psi (\epsilon d)\right)\notag
\end{eqnarray}
But, $c_k$ commutes with $x({\bf n})$. Hence, 
\begin{equation} \label{eqn1.4}
\left[ x({\bf n}), \, \varinjlim_{\epsilon \to 0}
L_{0}^{(k)}(\epsilon)\right] = \left[ x({\bf n}), \frac{1}{2}\sum_j\sum_{d\in \mathbbm{Z}} : e_j (-d) \cdot e^j (d):\right].
\end{equation}
Combining the equations \eqref{eqn1.2} and \eqref{eqn1.4}, we get the theorem. $\Square$
\end{proof}

\begin{remark} If ${\bf n}' = {\bf 0}$ in the above theorem, then $\left[x({\bf n}),L_0^{(k)}\right]$ as an operator on $V$ is given by the classical formula given in  \cite[Proposition 10.1]{KRR}.
\end{remark}


\begin{thebibliography}{99}

\bibitem[AJS]{AJS} H. Andersen, J. Jantzen and W. Soergel: Representations of quantum groups at a $p$-th root of unity and of semisimple groups in characteristic $p$ : independence of $p$, {\it Ast\'erisque}  {\bf 220}, 1-323 (1994).

\bibitem[F]{F} P. Fiebig, An upper bound on the exceptional characteristics for Lusztig's character formula, {\it Journal f\"ur die reine und angewandte Mathematik (Crelles Journal)} {\bf 673}, 1--31 (2012). 

\bibitem[GW]{GW}
  R. Goodman and N. Wallach,
  \emph{Symmetry,  Representations, and Invariants}.
  Graduate Texts in Mathematics Vol. \textbf{255},
  Springer (2009).


\bibitem[KRR]{KRR} V. Kac, A. Raina and N. Rezhkouskaya: {\it Bombay Lectures on highest weight representation of infinite dimensional Lie algebra (second edition)}, World Scientific (2013).

\bibitem[KT]{KT} M. Kashiwara and Tanisaki: Kazhdan-Lusztig conjecture for affine Lie algebras with negative level, {\it Duke Math. J.}  {\bf 77}, 21-- 62 (1995). 

\bibitem[KL0]{KL0} D. Kazhdan and G. Lusztig:  Affine Lie algebras and quantum groups, {\it International Math. Res. Notices}  {\bf 2}, 21--29 (1991).

\bibitem[KL1]{KL1} D. Kazhdan and G. Lusztig:  Tensor structures arising from affine Lie algebras I, {\it J. Amer. Math. Soc.}
{\bf 6}, 905--947 (1993).

\bibitem[KL2]{KL2} D. Kazhdan and G. Lusztig:  Tensor structures arising from affine Lie algebras II, {\it J. Amer. Math. Soc.}
{\bf 6}, 949--1011 (1993).

\bibitem[KL3]{KL3} D. Kazhdan and G. Lusztig:  Tensor structures arising from affine Lie algebras III, {\it J. Amer. Math. Soc.}
{\bf 7}, 335--381 (1994).

\bibitem[KL4]{KL2} D. Kazhdan and G. Lusztig:  Tensor structures arising from affine Lie algebras IV, {\it J. Amer. Math. Soc.}
{\bf 7}, 383--453 (1994).

\bibitem[K]{K}
  S. Kumar,
  \emph{Kac-Moody Groups, their Flag Varieties and  Representation Theory}.
  Progress in Mathematics Vol. \textbf{204}, Birkh\"auser   (2002).

\bibitem[Lu1]{Lu1} G. Lusztig:
Some problems in the representation theory of finite Chevalley groups, {\it Proc.
Symp. Pure Math.}  {\bf 37}, Amer. Math. Soc., 313--317 (1980).


\bibitem[Lu2]{Lu2} G. Lusztig:
On the character of certain irreducible modular representations, {\it Representation Theory, AMS} 
{\bf 19},  3--8 (2015).


\bibitem[Lu3]{Lu3} G. Lusztig: Comments on my papers, paper [229], arXiv:1707.09368 (v7) (2021).


\end{thebibliography}
\end{document}